\documentclass[12pt]{amsart}
\newtheorem{thm}{Theorem}[section]
\newtheorem{prop}[thm]{Proposition}
\newtheorem{lemma}[thm]{Lemma}
\newtheorem{coro}[thm]{Corollary}

\newtheorem{hypo}[thm]{Hypothesis {\bf H.}\hspace*{-0.6ex}}
\newtheorem{rem}[thm]{Remark}
\newcommand{\R}{{\mathbb R}}

\newcommand{\Z}{{\mathbb Z}}
\newcommand{\C}{{\mathbb C}}

\newcommand{\nn}{\nonumber}
\newcommand{\bea}{\begin{eqnarray}}
\newcommand{\eea}{\end{eqnarray}}
\newcommand{\ba}{\begin{array}}
\newcommand{\ea}{\end{array}}

\newcommand{\ol}{\overline}

\newcommand{\bth}{\begin{thm}}
\newcommand{\eth}{\end{thm}}
\newcommand{\bl}{\begin{lemma}}
\newcommand{\el}{\end{lemma}}
\newcommand{\bk}{\begin{coro}}
\newcommand{\ek}{\end{coro}}
\newcommand{\bh}{\begin{hypo}}
\newcommand{\eh}{\end{hypo}}
\newcommand{\br}{\begin{rem}}
\newcommand{\er}{\end{rem}}
\newcommand{\bpf}{\begin{proof}}
\newcommand{\epf}{\end{proof}}

\begin{document}

\title{Wiener Amalgam Spaces for the Fundamental Identity of Gabor Analysis}

\author{Hans G. Feichtinger}
\address{Fakult\"at f\"ur Mathematik\\
Nordbergstrasse 15\\ 1090 Wien\\ Austria}
\email{Hans.Feichtinger@univie.ac.at}

\author{Franz Luef}   
\address{Fakult\"at f\"ur Mathematik\\
Nordbergstrasse 15\\ 1090 Wien\\ Austria\\}
\address{Max Planck Institut f\"ur Mathematik,
Vivatsgasse 7\\
53111 Bonn\\
Germany} \email{Franz.Luef@univie.ac.at}

\maketitle \pagestyle{myheadings} \markboth{\Small Franz Luef and Hans G.
Feichtinger}{\Small Wiener Amalgam Spaces for the Fundamental Identity of
Gabor Analysis}

\begin{abstract}
In the last decade it has become clear that one of the central themes within
Gabor analysis (with respect to general time-frequency lattices) is a duality
theory for Gabor frames, including the Wexler-Raz biorthogonality condition,
the Ron-Shen's duality principle or Janssen's representation of a Gabor frame
operator. All these results are closely connected with the so-called {\it
Fundamental Identity of Gabor Analysis}, which we derive from an application
of Poisson's summation formula for the symplectic Fourier transform. The new
aspect of this presentation is the description of range of  the validity of
this Fundamental Identity of Gabor Analysis using Wiener amalgam spaces and
Feichtinger's algebra $S_0(\R^d)$. Our approach is inspired by Rieffel's use
of the Fundamental Identity of Gabor Analysis in the study of operator
algebras generated by time-frequency shifts along a lattice, which was later
independently rediscovered by Tolmieri/Orr, Janssen, and Daubechies et al.,
and Feichtinger/Kozek at various levels of generality, in the context of
Gabor analysis.
\end{abstract}

\section{Introduction}

Since the work of Wexler/Raz \cite{WR90} on the structure of the set of dual
atoms for a Gabor frame many researchers have benefited from their insight
that some properties of a Gabor frame have a better description with respect
to the adjoint lattice. We only mention the duality principle of Ron/Shen
\cite{RS93,RS97}, Janssen's representation of the Gabor frame operator
\cite{Jan95} and the investigations of Daubechies/H. Landau/Z. Landau
\cite{DLL95}, which have obtained similar results on the structure of Gabor
frames independently by completely different methods around the year 1995.
All their methods have in common an implicit use of the adjoint lattice for
separable lattices. In \cite{FK98}, Feichtinger and Kozek gained a thorough
understanding of the adjoint lattice for Gabor systems with respect to a
lattice $\Lambda$ in $G\times\widehat G$, $G$ an elementary locally compact
abelian group. Furthermore, \cite{FK98} makes use of the symplectic Fourier
transform in this context for the first time. Another ingredient in all
issues concerning investigations of duality principles for Gabor systems is
an identity about samples of the product of two short-time Fourier transforms
for a lattice $\Lambda$ and its adjoint lattice $\Lambda^0$, see Section
\ref{basicgab}. In Gabor analysis Tolimieri/Orr have pointed out the
relevance of this identity for the study of Gabor systems with atoms in the
Schwartz class \cite{TO95,TO92}. In \cite{Jan95}, Janssen generalized the results
of Tolimieri/Orr and called this identity the {\bf F}undamental {\bf
I}dentity of {\bf G}abor {\bf A}nalysis, since he derived FIGA as a
consequence of a representation of the Gabor frame operator of fundamental
importance in all duality results of Gabor analysis, (which is nowadays
called the Janssen representation of a Gabor frame operator). But, the
results of Tolimieri/Orr and Janssen on the FIGA had been obtained by Rieffel
in his construction of equivalence bimodules for the $C^*$-algebra $C^*(D)$
generated by time-frequency shifts of a closed subgroup $D$ of
$G\times\widehat G$, for a locally compact abelian group $G$, and the
$C^*$-algebra $C^*(D^0)$ generated by time-frequency shifts generated of the
adjoint group $D^0$ in $1988$, \cite{Rief88}. Furthermore, Rieffel gave a
description of the adjoint lattice, which was later rediscovered by
Feichtinger/Kozek, and he used the Poisson summation formula for the
symplectic Fourier transform to get FIGA for functions in the Schwartz-Bruhat
space ${\mathcal S}(G)$. Therefore, Tolimieri/Orr's discussion of the FIGA is
just a special case of Rieffel's general result. In addition Rieffel had
implicitly described Janssen's representation of a Gabor frame operator in
his discussion of $C^*(D^0)$-valued inner products  \cite{Rief88}.
\par
Our discussion of FIGA follows Rieffel's discussion. Therefore, we apply the
Poisson summation formula for the symplectic Fourier transform to a product
$V_{g_1}f_1\cdot\ol{V_{g_2}f_2}$ of short-time Fourier transforms of
functions resp.\ distributions $f_1,f_2,g_1,g_2$ in suitable modulation
spaces, see Section \ref{Modspaces} for the definition of modulation spaces.
In our proofs we need some local properties of STFT $V_{g_1}f_1$ and
$V_{g_1}f_1$ for $f_1,f_2,g_1,g_2$, which are naturally expressed by
membership in some Wiener amalgam spaces, see Section \ref{Modspaces}. Our
strategy relies heavily on the fact that Feichtinger's algebra $M^1(\R^{2d})$
is the biggest time-frequency homogenous Banach space, where the Poisson
summation formula holds pointwise (introduced as $S_0(\R^d)$ in
\cite{Fei81}). We therefore look for sufficient conditions such that
$V_{g_1}f_1\cdot\ol{V_{g_2}f_2}$ is in $M^1(\R^{2d})$.
\par
In Section \ref{basicgab} we introduce the reader to some well-known facts of
time-frequency analysis, which we will use later. In Section \ref{Modspaces}
we give a short discussion of modulation spaces and Wiener amalgam spaces and
present some of their properties. In Section \ref{figa} we discuss the notion
of weakly dual pairs in Gabor analysis and their connection to FIGA. Then we
prove FIGA with the help of Poisson summation for the symplectic Fourier
transform. Furthermore, we use Wiener amalgam spaces to describe the local
behaviour of the short-time Fourier transform. In Section \ref{S:Wexler-Raz} we
briefly indicate some consequences of our main result for Gabor frames.

\section{Basics of Gabor Analysis}\label{basicgab}

In this section we recall some well-known facts of Gabor analysis, e.g., the
short-time Fourier transform and some of its properties. Our representation
owes much to Gr\"ochenig's presentation in his survey of time-frequency
analysis \cite{GrBook}.
\par
In Gabor analysis the basic objects are time-frequency shifts. More
concretely, for $f\in L^2(\R^d)$  we define the following operators on
$L^2(\R^d)$:
\begin{enumerate}
  \item the {\it translation} operator by
\begin{equation*}
  T_xf(t)=f(t-x),\hspace{2.7cm}x\in\R^d,
\end{equation*}
  \item  the {\it modulation} operator by
\begin{equation*}
  M_{\omega}f(t)=e^{2\pi it\cdot\omega}f(t),\hspace{2.4cm}\omega\in\R^d,
\end{equation*}
  \item {\it time-frequency shifts} by
\begin{equation*}
  \pi(x,\omega)f=M_{\omega}T_xf=e^{2\pi i\omega t}f(t-x),\hspace{0.5cm}(x,\omega)\in\R^{2d}.
\end{equation*}
\end{enumerate}

The time-frequency shifts $(x,\omega,\tau)\mapsto\tau M_{\omega}T_x$ for
$(x,\omega)\in\R^{2d}$ and $\tau\in\C$ with $|\tau|=1$ define the
Schr\"odinger representation of the Heisenberg group, consequently the
time-frequency shifts $\pi(x,\omega)$ for $(x,\omega)\in\R^{2d}$ are a
projective representation of the time-frequency plane $\R^d\times{\widehat
\R}^d$. More concretely, time-frequency shifts satisfy the following
composition law:
\begin{equation}\label{TFcomp}
  \pi(x,\omega)\pi(y,\eta)=e^{-2\pi
  ix\cdot\eta}\pi(x+y,\omega+\eta),
\end{equation}
for $(x,\omega),(y,\eta)$ in the time-frequency plane $\R^d\times{\widehat
\R}^d$. The noncommutativity of the time-frequency shifts leads naturally to
the notion of the adjoint of a set of time-frequency shifts. Namely, let
$\Lambda$ be a subset of $\R^d\times{\widehat\R}^d$. Then, the {\it adjoint
lattice} $\Lambda^0$ of $\Lambda$ is defined as the set of all time-frequency
shifts in the time-frequency plane which commute with all time-frequency
shifts $\{\pi(\lambda):\lambda\in\Lambda\}$, i.e.,
\begin{equation}\label{adjointlattice}
  \Lambda^0:=\{\lambda^0\in\R^{2d}:\pi(\lambda)\pi(\lambda^0)=\pi(\lambda^0)\pi(\lambda)
  ~~~\text{for all}~~~\lambda\in\Lambda\}.
\end{equation}
\par
We include another approach to the adjoint of a lattice $\Lambda$ in
$\R^{d}\times{\widehat\R}^d$, because it plays a central role in the study of
Gabor frames.
\par
First we rewrite the composition law \eqref{TFcomp} of time-frequency shifts
\begin{equation}\label{TFcomm}
\pi(x,\omega)\pi(y,\eta)=e^{-2\pi i(x\cdot\eta-\omega\cdot
y)}\pi(y,\eta)\pi(x,\omega),
\end{equation}
for $(x,\omega),(y,\eta)$ in $\R^d\times{\widehat \R}^d$. We denote the
phase-factor in \eqref{TFcomm} by $\rho(z,z')=e^{2\pi i\Omega(z,z')}$ with
$z=(x,\omega)$, $z'=(u,\eta)$ and $\Omega$ denotes the standard symplectic
form on $\R^{2d}$, i.e. $\Omega(z,z')=x\cdot\eta-\omega\cdot y$. An important
fact is that $\rho$ is a character of $\R^d\times\widehat{\R}^d$ and that
every character of $\R^d\times\widehat{\R}^d$ is of the form
\begin{equation}
  z\mapsto\rho(z,z')~~~\text{for some}~~~z'\in\R^d\times\widehat{\R}^d.
\end{equation}
This gives an isomorphism between $\R^d\times\widehat{\R}^d$ and its dual
group $\widehat{\R}^d\times\R^d$.
\par
Let $\Lambda$ be a lattice in $\R^d\times\widehat{\R}^d$ then every character
of $\Lambda$ extends to a character of $\R^d\times\widehat{\R}^d$ and
therefore every character of $\Lambda$ is of the form
\begin{equation}
  \lambda\mapsto\rho(\lambda,z'),~~~\lambda\in\Lambda,
\end{equation}
for some $z\in\R^d\times\widehat{\R}^d$, where $z'$ needs not to be unique.
The homomorphism from $\R^d\times\widehat{\R}^d$ to $\widehat{\Lambda}$ has
as kernel:  the {\it adjoint lattice}
\begin{equation}
  \Lambda^0=\{z\in\R^d\times\widehat{\R}^d~|~~\rho(\lambda,z)=1
  ~~\text{for all}~~\lambda\in\Lambda\}.
\end{equation}
Therefore, the adjoint set $\Lambda^0$ of a lattice $\Lambda$ has the
structure of a lattice. In our discussion of FIGA we will explore further
this line of reasoning.
\par
The representation coefficients of the Schr\"odinger representation are up to
some phase factors, equal to
\begin{equation}
   V_gf(x,\omega):=\langle
   f,\pi(x,\omega)g\rangle=\int_{\R^d}f(t)\ol{g(t-x)}e^{-2\pi i\omega\cdot
   t}dt.
 \end{equation}
In Gabor analysis the representation coefficients are called the {\it
short-time Fourier transform} (STFT) of $f\in{\mathcal S}(\R^d)$ with respect
to a non-zero window $g$ in Schwartz's space of testfunctions ${\mathcal
S}(\R^d)$. For functions $f$ with good time-frequency concentration, e.g.
Schwartz functions, the STFT can be interpreted as a measure for the
amplitude of the frequency band near $\omega$ at time $x$. The properties of
STFT depend crucially on the window function $g$.
\par
In our model, the time-frequency concentration of a signal is invariant under
shifts in time and frequency, which is usually referred to the {\sl
covariance property} of a time-frequency representation. In harmonic
analysis, a function $f$ on $\R^d$ has a description in time and in
frequency. A time-frequency representation of a function $f$ encodes its
properties simultaneously in time and frequency, e.g. the STFT. The following
lemma expresses elementary properties of the STFT.
\begin{lemma}\label{basicSTFT}
 Let $f,g\in L^2(\R^d)$ and $(u,\eta)\in\R^d\times\widehat{\R}^d$. Then
  \begin{enumerate}
    \item {\it Covariance Property} of the STFT
          \begin{equation*}\label{covarianceSTFT}
            V_g(\pi(u,\eta)f)g(x,\omega)=e^{2\pi i
           u\cdot(\omega-\eta)}V_gf(x-u,\omega-\eta).
          \end{equation*}
    \item {\it Basic Identity of Time-Frequency Analysis}\label{basicident}
          \begin{equation*}\label{FI-STFT}
            V_gf(x,\omega)=e^{-2\pi ix\cdot\omega}V_{\hat{g}}\hat{f}(\omega,-x).
          \end{equation*}
  \end{enumerate}
\end{lemma}
In our proof of FIGA, we will use another basic identity of time-frequency
analysis: {\it Moyal's formula}.
\begin{lemma}[Moyal's Formula]\label{L:Moyal}
\noindent
  Let $f_1,f_2,g_1,g_2\in L^2(\R^d)$ then $V_{g_1}f_1$ and $V_{g_2}f_2$ are in $L^2(\R^{2d})$
  and the following identity holds:
  \begin{equation}\label{moyal}
    \langle V_{g_1}f_1,V_{g_2}f_2\rangle_{L^2(\R^{2d})}=\langle f_1,f_2\rangle
    \ol{\langle g_1,g_2\rangle}.
  \end{equation}
\end{lemma}
As a consequence, we get for $g\in L^2(\R^d)$ with $\|g\|_2=1$ that
\begin{equation*}
  \|V_gf\|_{L^2(\R^{2d})}=\|f\|_2,
\end{equation*}
for all $f\in L^2(\R^d)$, i.e., the STFT is an isometry from $L^2(\R^d)$ to
$L^2(\R^{2d})$.
\par
In time-frequency analysis we deal with function spaces which are invariant
under time-frequency shifts. In the last years {\it modulation spaces} have
turned out to be the correct class of Banach spaces for time-frequency
analysis, \cite{FK98,FZ98,GrBook,CG03}.

\section{Function Spaces for Time-Frequency Analysis}\label{Modspaces}
In the following we recall some well-known facts about modulation spaces and
Wiener amalgam spaces. Our treatment of this notions is largely based on the
excellent survey of time-frequency analysis by Gr\"ochenig, \cite{GrBook}.
\subsection{Modulation spaces}

In $1983$ Feichtinger introduced a class of Banach spaces (see
\cite{Fei83,Fei8302}), which allow a measurement of the time-frequency
concentration of a function or distribution $f$ on $\R^d$, the so called {\sl
modulation spaces}. We choose the STFT $V_gf$ of $f$ with respect to a window
$g$ with a good time-frequency concentration and as a measure we take the
norm of a function space which is(isometrically) invariant under translations
in the time-frequency plane $\R^d\times{\widehat \R}^d$. For our
investigations we restrict our study to weighted mixed-norm spaces
$L^{p,q}_m$ on $\R^{2d}$, \cite{Fei83}. But for the translation invariance of
$L^{p,q}_m$ Feichtinger showed that the weight $m$ has to be a moderate
weight on $\R^{2d}$ with respect to a positive and rotational symmetric
submultiplicative weight $v$ on $\R^{2d}$, i.e $m(z_1+z_2)\le Cm(z_1)v(z_2)$
for $z_1,z_2\in\R^{2d}$. Now for $1\le p,q\le\infty$ we define a function or
tempered distribution $f$ to be an element of the {\it modulation space}
$M^{p,q}_m(\R^d)$ if for a fixed $g$ in Schwartz space $\mathcal{S}(\R^d)$
the norm
\begin{equation*}
  \|f\|_{M^{p,q}_m}:=\|V_gf\|_{L^{p,q}_m}=
  \Big(\int_{\R^d}\Big(\int_{\R^d}|V_gf(x,\omega)|^p
  m(x,\omega)^pdx\Big)^{q/p}d\omega\Big)^{1/q}
\end{equation*}
is finite. Then $M^{p,q}_m(\R^d)$ is a Banach space whose definition is
independent of the choice of the window $g$. We always measure the
$M^{p,q}_m$-norm with a fixed non-zero window $g\in{\mathcal S}(\R^d)$ and
that for any non-zero $g_1\in M^1_v(\R^d)$ the norm equivalence of
$\|f\|_{M^{p,q}_m}$ with $\|V_{g_1}f\|_{L^{p,q}_m}$ holds.

\par
One reason for the usefulness of modulation spaces is that many well-known
function spaces can be identified with modulation spaces for certain weights:
\begin{enumerate}
  \item $M^{2,2}(\R^d)=L^2(\R^d)$.
  \item $M^1(\R^d)$ is {\it Feichtinger's algebra}, which is sometimes denoted by $S_0(\R^d)$.
  \item If $m(x,\omega)=(1+x^2)^{s/2}$ then $M^{2,2}_m=L^2_s=
  \{f\in{\mathcal
  S}'(\R^d):(\int_{\R^d}|f(x)|^2(1+x^2)^{s}dx)^{1/2}<\infty\}$ is a weighted $L^2$-space.
  \item If $m(x,\omega)=(1+\omega^2)^{s/2}$ then $M^{2,2}_m=H_s=
\{f\in{\mathcal S}'(\R^d):(\int_{\R^d}|\hat
f(\omega)|^2(1+\omega^2)^{s}d\omega)^{1/2}<\infty\}$ is a {\it Sobolev
space}.
  \item If $m(x,\omega)=(1+x^2+\omega^2)^{s/2}$ then $M^{2,2}_m=Q_s=L^2_s\cap H^s$,
  where $Q_s$ is the {\it Shubin class}, see \cite{Shu01}.
\end{enumerate}
Modulation spaces inherit many properties from the mixed norm spaces, e.g.
duality. In the following theorem we state some of their properties, that are
of interest in the later discussion.
\begin{thm}\label{ModspacesThm}
  Let $1\le p,q < \infty$ and $m$ a $v$-moderate weight on
  $\R^{2d}$.
  \begin{enumerate}
    \item The dual space of $M^{p,q}_m(\R^d)$ is $M^{p',q'}_{1/m}(\R^d)$ with
    $1/p+1/p'=1$ and $1/q+1/q'=1$ and the duality is given by
    \begin{equation*}
      \langle
      f,h\rangle=\iint_{\R^{2d}}V_{g}f(x,\omega)\ol{V_gh(x,\omega)}dxd\omega,
    \end{equation*}
    for $f\in M^{p,q}_m(\R^d)$ and $h\in M^{p',q'}_{1/m}(\R^d)$.
    \item \label{TF-invar} $M^{p,q}_m(\R^d)$ is invariant under time-frequency shifts:
\begin{equation*}
    \|\pi(u,\eta)f\|_{M^{p,q}_m}\le
    Cv(u,\eta)\|f\|_{M^{p,q}_m} \quad \mbox{for}  \quad  (u,\eta)\in\R^{2d}.
\end{equation*}
    \item If $p=q$ and $m(\omega,-x)\le Cm(x,\omega)$ then
    $M^{p,p}_m(\R^d)$ is invariant under Fourier transform.
 \end{enumerate}
\end{thm}
\begin{proof}
All these statements are well-known and the interested reader may find a
proof of statement $(1)$ in Chapter 11 of \cite{GrBook}. We only give the
arguments for statements $(2)$ and $(3)$, because they provide the reader
with some insight about our choice of weights.
\begin{enumerate}
  \item[(2)] The time-frequency invariance of $M^{p,q}_m(R^d)$ is a direct consequence of the definition of
  moderate weights and the Covariance Property of the STFT , see Lemma \ref{basicSTFT}. Let $z=(u,\eta)$ be
  a point of the time-frequency plane $\R^d\times\widehat{\R}^d$. Then the following holds:
  \begin{eqnarray*}
    \|\pi(u,\eta)f\|_{M^{p,q}_m}&=&\Big(\int_{\R^d}\Big(\int_{\R^d}|V_gf(x-u,\omega-\eta)|^p
  m(x,\omega)^pdx\Big)^{q/p}d\omega\Big)^{1/q}\\
                  &=&\Big(\int_{\R^d}\Big(\int_{\R^d}|V_gf(x,\omega)|^p
  m(x+u,\omega+\eta)^pdx\Big)^{q/p}d\omega\Big)^{1/q}\\
  &\le&C\Big(\int_{\R^d}\Big(\int_{\R^d}|V_gf(x,\omega)|^pv(u,\eta)^p
  m(x,\omega)^pdx\Big)^{q/p}d\omega\Big)^{1/q}\\
  &=&Cv(z)\|f\|_{M^{p,q}_m}.
 \end{eqnarray*}
\item[(3)] The key of the argument is an application of the basic identity of Gabor analysis, see Lemma \ref{basicSTFT}, to a
Fourier invariant window $g$ and the independence of the definition of
$M^{p,q}_m$ for $g\in{\mathcal S}(\R^d)$. For simplicity we choose $g$ to be
the standard Gaussian $g_0(x)=2^{-d/4}e^{-\pi x^2}$.
\begin{eqnarray*}
  \|\hat f\|_{M^{p,p}_m}&=&\Big(\iint_{\R^{2d}}|V_{g_0}f(x,\omega)|^p
  m(x,\omega)^pdxd\omega\Big)^{1/p}\\
    &\le&\Big(\iint_{\R^{2d}}|V_{\widehat{g_0}}\hat{f}(x,\omega)|^p
  m(x,\omega)^pdxd\omega\Big)^{1/p}\\
    &=&\Big(\iint_{\R^{2d}}|V_{g_0}f(-\omega,x)|^p
  m(x,\omega)^pdxd\omega\Big)^{1/p}\\
  &=&\Big(\iint_{\R^{2d}}|V_{g_0}f(x,\omega)|^p
  m(\omega,-x)^pdxd\omega\Big)^{1/p}\\
  &\le&C\|f\|_{M^{p,p}_m}.
\end{eqnarray*}

\end{enumerate}
\end{proof}
In the following Corollary, we state some of the properties of the modulation
space $M^{1,1}(\R^d)$. In harmonic analysis $M^{1,1}(\R^d)$ is the so-called
{\it Feichtinger algebra} and some authors use the notation $S_0(\R^d)$ to
indicate that Feichtinger's algebra is a Segal algebra, too. There is another
reason for this notation, because $S_0(\R^d)$ shares many properties with the
Schwartz space ${\mathcal S}(\R^d)$ of test functions, e.g., Feichtinger's
algebra is invariant under Fourier transform. In the rest of our paper we
will denote Feichtinger's algebra by $M^1(\R^d)$.
\begin{coro}
Feichtinger's algebra $M^1(\R^d)$ has the following properties:
\begin{enumerate}
 \item $M^1(\R^d)$ is a Banach algebra under pointwise multiplication.
 \item $M^1(\R^d)$ is a Banach algebra under convolution.
 \item $M^1(\R^d)$ is invariant under time-frequency shifts.
 \item $M^1(\R^d)$ is invariant under Fourier transform.
\end{enumerate}
\end{coro}
Before we present the proof we recall that the STFT can be written as a
convolution. Namely, let $g^*(x)=\ol{g(-x)}$ be the involution of $g\in
L^2(\R^d)$. Then, STFT of $f\in L^2(\R^d)$ has the following form
\begin{equation}\label{STFTconv}
   V_gf(x,\omega)=e^{-2\pi ix\cdot\omega}(f\ast M_{\omega}g^*)(x).
\end{equation}
For other formulations of the STFT and its relation to other time-frequency
representations, such as the Wigner distribution or the ambiguity function we
refer the reader to Gr\"ochening's book \cite{GrBook}.
\begin{proof}
  \begin{enumerate}
    \item By \eqref{STFTconv} the $M^1$-norm of $f$ is given by
      \begin{equation*}
         \|f\|_{M^1}=\int_{\R^d}\|f\ast M_{\omega}g^*\|_{L^1}d\omega.
      \end{equation*}
      Therefore, we get the following estimate for $h\ast f$, where $f\in M^1(\R^d)$ and $h\in L^1(\R^d)$:
      \begin{eqnarray*}
         \|h\ast f\|_{M^1}&=&\int_{\R^d}\|h\ast f\ast M_{\omega}g^*\|_{L^1}d\omega\\
     &\le&\int_{\R^d}\|h\|_{L^1}\|f\ast M_{\omega}g^*\|_{L^1}d\omega\\
     &=&\|h\|_{L^1}\|f\|_{M^1}.
      \end{eqnarray*}
      \item The statement follows from $(1)$ by applying Fourier transforms.
      \item The statement is a special case of our general result for modulation spaces, Theorem \ref{ModspacesThm}.
      \item The statement is again a special case of our general result for modulation spaces, Theorem \ref{ModspacesThm}.
    \end{enumerate}
\end{proof}
Despite the above stated properties, Feichtinger observed that $M^1(\R^d)$ is
the {\it minimal} time-frequency homogenous Banach space, \cite{Fei81}.
Another pleasant property of Feichtinger's algebra $M^1(\R^d)$ is that it is
the largest Banach space which allows an application of Poisson's summation
formula, \cite{Fei81}. Our main results about FIGA rely heavily on this fact.
We will discuss this topic further after the introduction of the symplectic
Fourier transform.
\par

\subsection{Wiener Amalgam Spaces}
Around $1980$ Feichtinger introduced a class of Banach spaces, which allow
measurement of local and of global properties of functions, see
\cite{Fei81a}. Feichtinger's work was motivated by some spaces Wiener had
used in his study of the Fourier transform, see \cite{Wie}. Nowadays those
spaces are called {\it Wiener amalgam spaces} and they are a generalization
of Fournier and Stewart's amalgam spaces, \cite{FS85}. Wiener amalgam spaces
have turned out to be very useful in harmonic analysis and time-frequency
analysis, e.g. \cite{Fei90,FZ98,GrBook}.
\par
More concretely, let $g$ be an element of the space of test functions
${\mathcal D}(\R^d)$, whose translates generate a partition of unity over
$\R^{2d}$, i.e., $\sum_{m\in\Z^{2d}}T_{m}g\equiv 1$. Let $V(\R^{2d})$ be a
translation invariant Banach space of functions (or distributions) over
$\R^{2d}$ with the property that $\mathcal{D}\cdot V\subset V$. Then the {\it
Wiener amalgam space} $W(X,L^{p,q}_m)$ with local component $X$ and global
component $L^{p,q}_m$ is defined as the space of all functions resp.\
distributions for which the norm
\begin{equation*}
  \|f\|_{W(X,L^{p,q}_m)}=\Big(\int_{\R^d}\big(\int_{\R^d}\|f\cdot T_{(z_1,z_2)}\ol{g}\|_X^p
  m(z_1,z_2)^pdz_1\big)^{q/p}dz_2\Big)^{1/q}
\end{equation*}
is finite. We note that different choices of $g\in{\mathcal D}$ give the same
space and yield equivalent norms \cite{Fei81a}.
\par
In \cite{Fei81a} also an extension of H\"older's inequality to Wiener amalgam
spaces is given: If  $X$ be a Banach algebra with respect to pointwise
multiplication, then for $\frac{1}{p}+\frac{1}{p'}=1$ and
$\frac{1}{q}+\frac{1}{q'}=1$ ones has
\begin{equation}\label{hoelder}
  \|f\|_{W(X,L^1)}\le\|f\|_{W(X,L^{p,q}_m)}\|f\|_{W(X,L^{p',q'}_{1/m})}\quad .
\end{equation}
\par
There are many characterizations of Feichtinger's algebra $M^1(\R^d)$. Later
we shall need a result of Feichtinger that $M^1(\R^d)=W({\mathcal F}L^1,L^1)$
with equivalent norms \cite{Fei81}. One of these norms
has a formulation by means of the STFT:
\begin{eqnarray*}
  \|f\|_{W({\mathcal F}L^1,L^1)}&=&\iint_{\R^{2d}}\|f\cdot T_{(z_1,z_2)}g\|_{{\mathcal F}L^1}
  dz_1dz_2\\
  &=&\iint_{\R^{2d}}\Big(\int_{\R^d}|(f\cdot T_{(z_1,z_2)}\ol{g})^{\widehat{ }}(\omega)|d\omega\Big)
  dz_1dz_2\\
  &=&\iint_{\R^{2d}}\Big(\int_{\R^d}|V_gf(z_1,z_2,\omega)|d\omega\Big)dz_1dz_2.
\end{eqnarray*}
Analogous expressions for the norm of Wiener amalgam spaces $W({\mathcal F}
L^1,L^{p,q}_m)$ with local component ${\mathcal F}L^1$ and global component
$L^{p,q}_m$ can be derived in terms of the STFT.

\section{The Fundamental Identity of Gabor Analysis}\label{figa}

In this section we first recall some elementary facts about Gabor frame
operators and the symplectic Fourier transform. These basic facts and a
result \cite{CG03} of Cordero/Gr\"ochenig on local properties for $V_gf$ with
window $g\in M_v^1(\R^d)$ and $f\in M^{p,q}_m(\R^d)$ will allow us to derive
our main result about the Fundamental Identity of Gabor Analysis.
\par
Let $\Lambda$ be a lattice in $\R^d\times{\widehat\R}^d$ and let $g\in
L^2(\R^d)$ be a Gabor atom then ${\mathcal
G}(g,\Lambda)=\{\pi(\lambda)g:\lambda\in\Lambda\}$ is a {\it Gabor system}.
Since the Balian-Low principle tells us that it is not possible to
construct (this is in contrast to the situation with wavelets) an
orthonormal basis for $L^2(\R^d)$ of this form, starting e.g.\ from a
Schwartz function $g$,
interest in Gabor frames arose. A milestone was the paper
by Daubechies, Grossmann  and Meyer,
\cite{DGM86}, where the ``painless use'' of (tight) Gabor frames was
suggested. \par
In our case the Gabor frame operator has the following form:
\begin{equation*}
  S_{g,\Lambda}f=\sum_{\lambda\in\Lambda}\quad \langle f,\pi(\lambda)g\rangle\pi(\lambda)g,\quad  f  \in L^2(\R^d).
\end{equation*}
A Gabor system ${\mathcal G}(g,\Lambda)$ is called a {\it Gabor frame} if the
Gabor frame operator $S_{g,\Lambda}$ is invertible, i.e., if there exist some
finite, positive real numbers $A,B$ such that
\begin{equation*}
  A\cdot{\text I}\le S_{g,\Lambda}\le B\cdot{\text I}
\end{equation*}
or equivalently,
\begin{equation*}
  A\|f\|^2\le\sum_{\lambda\in\Lambda}|\langle f,\pi(\lambda)g\rangle|^2\le B\|f\|^2,
\end{equation*}
for all $f$ in $L^2(\R^d)$. Gabor frames ${\mathcal G}(g,\Lambda)$ allow the
following reconstruction formulas
\begin{eqnarray}\label{reconstruct}
    f&=&(S_{g,\Lambda})^{-1}S_{g,\Lambda}f=\sum_{\lambda\in\Lambda}\langle
  f,\pi(\lambda)g\rangle\pi(\lambda)(S_{g,\Lambda})^{-1}g\\
  f&=&S_{g,\Lambda}(S_{g,\Lambda})^{-1}f=\sum_{\lambda\in\Lambda}\langle
  f,\pi(\lambda)(S_{g,\Lambda})^{-1}g\rangle\pi(\lambda)g.
\end{eqnarray}
Due to its appearance in the reconstruction formulas
$\gamma_0:=(S_{g,\Lambda})^{-1}g$ is called
the (canonical) 
{\it dual Gabor atom}. Note that the non-orthogonality of the time-frequency
shifts yields that the coefficients in the reconstruction formula
\eqref{reconstruct} are not unique and therefore there are other {\it dual
atoms} $\gamma\in L^2(\R^d)$ with
$S_{g,\gamma,\Lambda}:=\sum_{\lambda\in\Lambda}\langle
f,\pi(\lambda)\gamma\rangle\pi(\lambda)g={\text I}$. Some authors call
$(g,\gamma)$ a {\it dual pair} of Gabor atoms if $S_{g,\gamma,\Lambda}={\text
I}$.
\par
If one considers Gabor systems ${\mathcal G}(\Lambda,g)$ beyond
$L^2$-setting, then the operator identity $S_{g,\gamma,\Lambda}={\text I}$
has to be interpreted in weak sense. This approach was suggested by
Feichtinger and Zimmermann in \cite{FZ98}. They called two elements
$g,\gamma\in L^2(\R^d)$  
a {\it weakly dual pair} with respect to
$\Lambda$, if
\begin{equation}
  \langle f,h\rangle=\sum_{\lambda\in\Lambda}\langle f,\pi(\lambda)\gamma\rangle\langle\pi(\lambda)g,h\rangle
\end{equation}
holds. One can show that absolute convergence of the series on the right,
for all $f,h\in M^1(\R^d)$.
Although the absolute convergence of the series seems very restrictive at first
sight
also for the case that  $g\in M^1(\R^d)$ and $\gamma\in M^{\infty}(\R^d)$. By the
symmetry of the definition in $g$ and $\gamma$ we get the same result for
$\gamma\in M^1(\R^d)$ and $g\in M^{\infty}(\R^d)$. Therefore, in the sequel
we will only state our results for one setting. We also mention without proof
that under the above assumption on the pair $(g,\gamma)$ the operator mapping
the pair $(f,h)$ to $\sum_{\lambda\in\Lambda}\langle
f,\pi(\lambda)\gamma\rangle\langle\pi(\lambda)g,h\rangle$ is continuous on
$M^1(\R^d)\times M^1(\R^d)$ and that the corresponding Gabor frame operator
$S_{g,\gamma,\Lambda}$ maps $M^1(\R^d)$ into $M^{\infty}(\R^d)$, \cite{FZ98}.
Another direct consequence of the definiton is that a pair $(g,\gamma)$ is
weakly dual if and only if $S_{g,\gamma,\Lambda}={\text I}$. For further
properties we refer the interested reader to \cite{FZ98}.
\par
The preceding discussion leads in a natural way to the study of the product
of two STFT's $V_{\gamma}f\cdot\ol{V_gh}$ restricted to a lattice $\Lambda$.
In Gabor analysis Tolmieri and Orr realized (in the one-dimensional case and
for a product lattice $\alpha\Z\times\beta\Z$, $\alpha,\beta>0$) that such
sums should be evaluated with the help of Poisson's summation formula. But in
his research on Morita equivalence of noncommuative tori \cite{Rief88},
Rieffel had used this identity for functions $f,g,h,\gamma$ in
Schwartz-Bruhat space ${\mathcal S}(G)$ for an elementary locally compact
abelian group $G$ and restrictions to a closed subgroup $D$ of
$G\times\widehat G$ in $1988$! Only recently one of us has realized the
connection between Rieffel's results and Gabor analysis \cite{Lu05}.
Therefore the research in Gabor anlaysis has been undertaken independently of
Rieffel's work, despite its great relevance for Gabor anlaysis. We will
explore this further in a subsequent paper.
\par
In the following we define the symplectic Fourier transform and some of its
basic properties, which was implicitely used by Rieffel in his derivation of
the Fundamental Identity of Gabor Analysis. Due to the work of Feichtinger
and Kozek \cite{FK98}  we also have gained the insight that the symplectic
Fourier transform might be of some relevance in this context.
\par
We continue our investigations of the character $\rho$ of the commutation
relation \eqref{TFcomm}. The antisymmetry of $\Omega$ implies that $\rho$ is
a skew-bicharacter of $\R^d\times\widehat{\R}^d$. Nevertheless, $\rho$ gives
a Fourier transform ${\widehat F}^s$ on the time-frequency plane
$\R^d\times{\widehat\R}^d$
\begin{eqnarray*}\nonumber
{\widehat F}^s(z)&=&\iint_{\R^{2d}}\rho(z,z')F(z')dz'\\
                 &=&\iint_{\R^{2d}}e^{2\pi i\Omega(z,z')}F(z')dz'\\
                 &=&\iint_{\R^{2d}}e^{2\pi i(y\cdot\omega - x\cdot\eta)}F(y,\eta)dyd\eta,
\end{eqnarray*}
for $z=(x,\omega)$ and $z'=(y,\eta)$ in $\R^d\times{\widehat\R}^d$. We call
${\widehat F}^s$ the {\it symplectic Fourier transform} of a function $F$ in
$L^2(\R^d\times\widehat{\R}^d)$, because it is induced by the symplectic form
$\Omega$ of $\R^{2d}$.
\par
The Poisson summation formula is one of the most powerful tools in harmonic
analysis. In our derivation of FIGA we need a Poisson summation formula for
the symplectic Fourier transform, which relates values of a function $F$ on a
lattice $\Lambda$ in the time-frequency plane with the samples of its
symplectic Fourier transform on the adjoint lattice $\Lambda^0$.
\par
In the following theorem we give some properties of the symplectic Fourier
transform.
\begin{thm}
Let $M^1(\R^{2d})$ be Feichtinger's algebra over the time-frequency plane
$\R^d\times{\widehat\R}^d$.
  \begin{enumerate}
    \item The symplectic Fourier transform is selfinverse on $L^2(\R^{2d})$.
    \item $M^1(\R^{2d})$ is invariant under the symplectic Fourier transform.
  \end{enumerate}
\end{thm}
\begin{proof}
\noindent
  \begin{enumerate}
     \item The statement is a consequence of the anti-symmetry of the symplectic form $\Omega$.
     \item We make the observation that $\rho$ arises from $e^{2\pi i(x\cdot\omega+y\cdot\eta)}$ by a rotation of
           $\pi/2$, i.e. the symplectic Fourier transform is a rotated version of the Fourier transform on $\R^{2d}$.
       Therefore the result follows from the main properties of Feichtinger's algebra, see Theorem \ref{ModspacesThm}.
  \end{enumerate}
\end{proof}

Traditionally a harmless use of the Poisson summation formula is only known
for Schwartz functions. That for Feichtinger's algebra Poisson summation
holds pointwise and with absolute convergence is quite unexpected. We only
remind you on the work of Katznelson and of Kahane et al., where they give
striking examples for the failure of Poisson's summation formula,
\cite{KL94,Kat67}. In \cite{Gr96} Gr\"ochenig pointed out the relevance of
Poisson's summation formula for $M^1(\R^d)$ in his study of uncertainty
principles and embeddings of various function spaces into $M^1(\R^d)$.
\par
In the following theorem we state Poisson's summation formula for the
symplectic Fourier transform.
\begin{thm}\label{PoissonSum}[Poisson Summation]
Let $F\in M^1(\R^{2d})$ then
\begin{equation}
  \sum_{\lambda\in\Lambda}F(\lambda)=|\Lambda|^{-1}\sum_{\lambda^0\in\Lambda^0}
  {\widehat{F}}^s(\lambda^0)
\end{equation}
holds pointwise and with absolute convergence of both sums.
\end{thm}
In \cite{BP04} Benedetto/Pfander constructed a function $g\in L^2(\R)$ such
that $|V_gg|^2\notin M_1(\R^2)$ and therefore the symplectic Fourier
transform is not valid for $|V_gg|^2$.

 Before we present our results on the validity of the FIGA, we compute the symplectic Fourier transform
 of $V_gf$ for $f,g\in M^1(\R^{d})$.
 \begin{lemma}
 Let $f,g\in M^1(\R^d)$ then the following holds:
 \begin{enumerate}
    \item $V_gf\in M^1(\R^{2d})$.
    \item $\widehat{V_gf}^s(z)=f(x)\ol{\hat g(\omega)}e^{-2\pi ix\cdot\omega}$ for $z=(x,\omega)\in\R^{2d}$.
 \end{enumerate}
 \end{lemma}
\begin{proof}
  \begin{enumerate}
  \item   Follows from the functorial properties and the minimality of $M^1(\R^d)$, see \cite{Fei81}. Feichtinger and
   Kozek give a different proof in \cite{FK98}.
   \item Straightforward computation.
  \end{enumerate}
\end{proof}
We recall that {\it Rihaczek's distribution} for $f,g\in L^2(\R^d)$ is
defined as $R(f,g)(x,\omega)=f(x)\ol{\hat g(\omega)}e^{-2\pi ix\cdot\omega}$ for
$z=(x,\omega)\in\R^{2d}$. It is a very popular
time-frequency representation in engineering, \cite{HM01,Gr04}. An
application of Theorem \ref{PoissonSum} yields a generalization of a formula
of Kaiblinger \cite{Ka05} 
which also established a connection to the Rihaczek distribution.
\begin{prop}
\noindent
 Let $f,g$ be in $M^1(\R^d)$ and $\Lambda$ a
 lattice in $\R^d\times{\widehat\R}^d$. Then the
 following relation holds:
 \begin{equation*}
   \sum_{\lambda\in\Lambda}V_gf(\lambda)=|\Lambda|^{-1}\sum_{\lambda^0\in\Lambda^0}
   R(f,g)(\lambda^0),
 \end{equation*}
\end{prop}

\par
An application of Theorem \eqref{PoissonSum} to a product of two STFT's
combined with Lemma (4.3)(2) yields the FIGA.

\begin{thm}[Basic FIGA]
Assume that $f_1,f_2,g_1,g_2\in M^1(\R^d)$. Then
\begin{equation}\label{FIGA}
  \sum_{\lambda\in\Lambda}V_{g_1}f_1(\lambda) \cdot \overline{V_{g_2}f_2(\lambda)}=
  |\Lambda|^{-1}\sum_{\lambda^0\in\Lambda^0}    V_{g_1}g_2(\lambda^0) \cdot
  \overline{V_{f_1}f_2(\lambda^0)}
\end{equation}
\end{thm}
\begin{proof}
 Our argument is just the computation of the symplectic Fourier transform of
 $F=V_{g_1}f_1(\lambda) \cdot \overline{V_{g_2}f_2(\lambda)}$ and a use of Theorem \eqref{PoissonSum}.
 \begin{eqnarray}\nn
   \widehat{F}^s(Y)&=&\iint_{\R^d\times{\widehat{\R}}^d}
                      V_{g_1}f_1(X)\overline{V_{g_2}f_2(X)}\rho(Y,X)dY\\\nn
                   &=&\iint_{\R^d\times{\widehat{\R}}^d}
                   \langle\pi(Y)f_1,\pi(Y)\pi(X)g_1\rangle
                   \overline{\langle
                   f_2,\pi(X)g_2\rangle}\rho(Y,X)dY\\\nn
                   &=&\iint_{\R^d\times{\widehat{\R}}^d}
                   \langle\pi(Y)f_1,\pi(X)\pi(Y)g_1\rangle
                   \overline{\langle
                   f_2,\pi(X)g_2\rangle}\rho(X,Y)dY\\\nn
                   &=&\langle f_1,\pi(Y)f_2\rangle\overline{\langle
                   g_1,\pi(Y)g_2\rangle},
 \end{eqnarray}
 where in the last step we used Moyal's formula \eqref{moyal}.
 \end{proof}
 The Fourier invariance of $M^1(\R^d)$ and the basic identity of Gabor analysis \eqref{basicident} yield the following reformulation of
 \eqref{FIGA}:

\begin{equation*}
  \sum_{\lambda\in{\mathcal J}\Lambda}V_{\hat{g_1}}\hat{f_1}(\lambda) \cdot \overline{V_{\hat{g_2}}\hat{f_2}(\lambda)}=
  |\Lambda|^{-1}\sum_{\lambda^0\in{\mathcal J}\Lambda^{\perp}} V_{g_1}g_2(\lambda^0) \cdot,
  \overline{V_{f_1}f_2(\lambda^0)}
\end{equation*}
  where ${\mathcal J}=\left( \begin{array}{cc}0 &{\text I}_d \\-{\text I}_d & 0 \end{array} \right)$ denotes a rotation by $\pi/2$ of the
  time-frequency plane $\R^d\times{\widehat \R}^d$. We have stated this reformulation of \eqref{FIGA} for two reasons:
  (1) It is another manifestation of the fact that $V_gf$ encodes information about $f$ and $\hat f$,(2) It shows that an application of
  the Fourier transform on the level of signals and windows corresponds to a rotation of the lattice $\Lambda$ in $\R^d\times{\widehat\R}^d$ and
  that the application of the symplectic Fourier transform to a time-frequency representation of $f$ yields a rotation of the dual lattice
  $\Lambda^{\perp}$ in ${\widehat\R}^d\times\R^d$.
\par
What properties of $f,g\in M_1(\R^d)$ imply that $V_gf\in M_1(\R^{2d})$?
Recall that $M_1(\R^d)=W({\mathcal F}L^1,L^1)$. This fact suggests that for
$f\in M^{p,q}_m(\R^d)$ and $g\in M^1(\R^d)$ the STFT $V_gf\in W({\mathcal
F}L^1,L^{p,q}_m)$. In \cite{CG03} Cordero and Gr\"ochenig have proved this
result in their investigations of localization operators. Before we state
their result on local properties of the STFT, we introduce a weighted version
of Feichtinger's algebra $M^1(\R^d)$.
\par
By assumption our weight $m$ is $v$-moderate for $v$ a submultiplicative
weight on $\R^{2d}$. The space $M^1_v(\R^d)$ is now the correct weighted
version of Feichtinger's algebra $M^1(\R^d)$.

\begin{prop}[Cordero/Gr\"ochenig]\label{localSTFT}
  Let $1\le p,q\le\infty$. If $f\in M^{p,q}_m(\R^d)$ and $g\in
  M^1_v(\R^d)$, then $V_gf\in W({\mathcal F}L^1,L^{p,q}_m)$ with
  \begin{equation}\label{localSTFTnorm}
    \|V_gf\|_{W({\mathcal F}L^1,L^{p,q}_m)}\le
    C\|f\|_{M^{p,q}_m}\|g\|_{M^1_v}.
  \end{equation}
 \end{prop}

In other words, the norms of $L^{p,q}_m$ and of $W({\mathcal
F}L^1,L^{p,q}_m)$ are equivalent on the range of $V_g$.
\par
We refer the reader to \cite{CG03} for a proof of Proposition
\ref{localSTFT}.
\par

The remaining part of this section we look for conditions on the quadruple
$f_1,f_2,g_1,g_2$, which imply that $V_{g_1}f_1 \cdot   \ol {V_{g_2}f_2}$,
belongs to Feichtinger's algebra $M^1(\R^d)$, because then a careless
application of Poisson formula summation is allowed.
\bigskip


The following theorem is our main result, which is an extension of the
validity of range of FIGA.

\begin{thm}\label{mainthm}[Main Result] There exists $C > 0$ independent of $p,q,v,m$ such that for
   $f_1\in M^{p,q}_m$, $f_2\in M^{p',q'}_{1/m}$ and $g_1,g_2\in M^1_v$, then
   the following holds:
 \begin{equation*}
  \sum_{\lambda\in\Lambda}V_{g_1}f_1(\lambda) \cdot \overline{V_{g_2}f_2(\lambda)}=
  |\Lambda|^{-1}\sum_{\lambda^0\in\Lambda^0}    V_{g_1}g_2(\lambda^0) \cdot
  \overline{V_{f_1}f_2(\lambda^0)}.
\end{equation*}
\end{thm}
\begin{proof}
  By Proposition \ref{localSTFT} we have that $V_{g_1}f_1\in W({\mathcal
  F}L^1,L^{p,q}_m)$ and $V_{g_2}f_2\in W({\mathcal
  F}L^1,L^{p',q'}_{1/m})$. Therefore an application of H\"older's
  inequality \eqref{hoelder} for Wiener amalgam spaces yields that $V_{g_1}f_1\cdot\ol{V_{g_2}f_2}\in  W({\mathcal
  F}L^1,L^1)$. The inequalities \eqref{hoelder} and \eqref{localSTFTnorm}
  imply the desired norm estimate:
  \begin{eqnarray*}
    \|V_{g_1}f_1\cdot\ol{V_{g_2}f_2}\|_{M^1}&\le&
    \|V_{g_1}f_1\cdot\ol{V_{g_2}f_2}\|_{W({\mathcal
  F}L^1,L^1)}\\
  &\le&C\|V_{g_1}f_1\|_{W({\mathcal
  F}L^1,L^{p,q}_m)}\|V_{g_2}f_2\|_{W({\mathcal
  F}L^1,L^{p',q'}_{1/m})}\\
  &\le&
  C\|g_1\|_{M^1_v}\|g_2\|_{M^1_v}\|f_1\|_{M^{p,q}_m}\|f_2\|_{M^{p',q'}_{1/m}}.
  \end{eqnarray*}
  Therefore our object of interest is in $M^1(\R^d)$ and an application of Poisson
  summation yields the desired result.
\end{proof}

As an application of Theorem \ref{mainthm} we derive the known results about
the validity of the FIGA. The first result was obtained by
Feichtinger/Zimmermann in their discussion of weakly dual pairs \cite{FZ98}.
\begin{coro}[Feichtinger-Zimmermann]
Let $g_1,g_2$ be in $M^1(\R^d)$. If $f_1\in M^1(\R^d)$ and $f_2\in
M^{\infty}(\R^d)$ or $f_1,f_2\in L^{2}(\R^d)$ then
   \begin{equation*}
  \sum_{\lambda\in\Lambda}V_{g_1}f_1(\lambda) \cdot \overline{V_{g_2}f_2(\lambda)}=
  |\Lambda|^{-1}\sum_{\lambda^0\in\Lambda^0}    V_{g_1}g_2(\lambda^0) \cdot
  \overline{V_{f_1}f_2(\lambda^0)}.
\end{equation*}
    holds.
\end{coro}
\begin{proof}
 The corollary covers the cases $f_1\in M^1(\R^d)$ and $f_2\in \big(M^1(\R^d)\big)'=M^{\infty}(\R^d)$ and $f_1,f_2\in M^{2,2}(\R^d)=L^2(\R^d)$.
 Therefore the proof is a direct consequence of Theorem \ref{mainthm}
\end{proof}
The second result covers the case of Tolimieri/Orr of the validity of the
FIGA for Schwartz funtions \cite{TO95}. The proof consists of the well-known
fact \cite{GrBook} that the modulation spaces $M^1_{v_s}(\R^d)$ for
$v_s(x,\omega)=(1+x^2+\omega^2)^{s/2}$ are the building blocks of the
Schwartz class ${\mathcal S}(\R^d)$, namely $${\mathcal
S}(\R^d)=\bigcap_{s\ge 0}M^1_{v_s}(\R^d).$$ By duality we get a description
of tempered distributions $${\mathcal S}'(\R^d)=\bigcup_{s\ge
0}M^{\infty}_{1/v_s}(\R^d).$$

\begin{coro}[Tolimieri-Orr]
Let $f_1,g_1,g_2$ be in ${\mathcal S}(\R^d)$ and $f_2\in{\mathcal S}'(\R^d)$
then we have the following identity:
 \begin{equation*}
  \sum_{\lambda\in\Lambda}V_{g_1}f_1(\lambda) \cdot \overline{V_{g_2}f_2(\lambda)}=
  |\Lambda|^{-1}\sum_{\lambda^0\in\Lambda^0}    V_{g_1}g_2(\lambda^0) \cdot
  \overline{V_{f_1}f_2(\lambda^0)}.
\end{equation*}

\end{coro}
\begin{proof}
The statement is true for every building block of ${\mathcal S}(\R^d)$ and of
${\mathcal S}'(\R^d)$, respectively. Therefore our statement is a direct
consequence of our main result Theorem \ref{mainthm}.
\end{proof}

We stop here our list of examples and leave it to the reader to choose a
pairing of his interest.

\section{Biorthogonality condition of Wexler-Raz}\label{S:Wexler-Raz}

In this section we present some consequences of our results on the FIGA,
especially its relation to Janssen's representation of Gabor frame operators
and the biorthogonality condition of Wexler-Raz.
\par
Many researchers have drawn deep consequences from Janssen's representation,
e.g. Gr\"ochenig/Leinert in their proof of the "irrational case"-conjecture
\cite{GL04}, Feichtinger/Kaiblinger in their work on the continuous
dependence of the dual atom for a Gabor atom in $M^1(\R^d)$ or ${\mathcal
S}(\R^d)$ \cite{FK04}.
\par
In \cite{Jan95} Janssen obtained his representation for Gabor frame operators
$S_{g,\Lambda}$ with $g\in{\mathcal S}(\R^d)$ and
$\Lambda=\alpha\Z\times\beta\Z^d$. The general form for arbitrary lattices
was obtained by Feichtinger in collaboration with Kozek and Zimmermann in
\cite{FK98,FZ98}.
\par
Let $\Lambda$ be a lattice in $\R^d\times\widehat{\R}^d$ and $g$ a Gabor atom
then the frame operator
\begin{equation*}
   S_{g,\Lambda}f=\sum_{\lambda\in\Lambda}\langle f,\pi(\lambda)g\rangle\pi(\lambda)g.
\end{equation*}
Janssen's insight consists on a formal level of the following observation

\begin{eqnarray*}
  \langle S_{g,\Lambda}f,h\rangle&=&\Big\langle\sum_{\lambda\in\Lambda}\langle f,\pi(\lambda)g\rangle\pi(\lambda)g,h\Big\rangle\\
                                 &\stackrel{(FIGA)}{=}&\Big\langle|\Lambda|^{-1}\sum_{\lambda^0\in\Lambda^0}\langle g,\pi(\lambda^0)g\rangle\pi(\lambda^0)f,h\Big\rangle\\
                 &=&\langle|\Lambda|^{-1} S_{f,\Lambda^0}g,h\Big\rangle.
\end{eqnarray*}
But $S_{f,\Lambda^0}g$ is a series of time-frequency shifts operators acting
on $f$. More concretely, the following representation of the Gabor frame
operator was obtained by Janssen in \cite{Jan95}
\begin{equation}\label{Janssenrep}
   S_{g,\Lambda}=|\Lambda|^{-1}\sum_{\lambda^0\in\Lambda^0}\langle g,\pi(\lambda^0)g\rangle\pi(\lambda^0).
\end{equation}
But the series on the right side of \eqref{Janssenrep} only defines a bounded
operator with the additional assumption that
\begin{equation}\label{condA}
 \text{(A)}\hspace{2.5cm}  \sum_{\lambda^0\in\Lambda^0}|\langle g,\pi(\lambda^0)g\rangle|<\infty.
\end{equation}
The last condition was introduced by Tolimieri/Orr in their discussion of
Gabor frames \cite{TO95}.
\par
The preceding observations led Janssen to consider operators of the form
\begin{equation*}
     A=\sum_{\lambda^0\in\Lambda^0}a(\lambda^0)\pi(\lambda^0)
\end{equation*}
for $\big(a(\lambda^0)\big)\in\ell^1(\Lambda^0)$. But in \cite{Rief88}
Rieffel used such operators to introduce on ${\mathcal S}(\R^d)$ a Hilbert
$C^*$-module structure for the $C^*$-algebra of all time-frequency shifts
generated by $\Lambda^0$. This fact is the reason for the relation between
Rieffel's work on Morita equivalence of noncommutative tori and Gabor
analysis \cite{Lu05}.
\par
We now extend the main results of Feichtinger/Zimmermann about weakly dual
pairs \cite{FZ98} to our setting. The notion of weakly dual pairs is the
proper concept for the interpretation of a Gabor frame operator
$S_{g,\gamma,\Lambda}$ in a weak sense.
\par
First we recall the biorthogonality condition of Wexler-Raz. In Section
\ref{figa} we have shown that the Gabor frame operator $S_{g,\Lambda}$ of a
Gabor system ${\mathcal G}(g,\Lambda)$ for a lattice $\Lambda\in\R^{2d}$
gives rise to a reconstruction formula
\begin{equation}
f=(S_{g,\Lambda})^{-1}S_{g,\Lambda}f=\sum_{\lambda\in\Lambda}\langle
  f,\pi(\lambda)g\rangle\pi(\lambda)\gamma_0
\end{equation}
for signals $f\in L^2(\R^d)$. We also mentioned the existence of other dual
functions $\gamma$, which give rise for reconstruction formulas. In
\cite{WR90} Wexler/Raz gave a characterization of all dual functions $\gamma$
for periodic discrete Gabor systems, which was the motivation for the work of
Janssen, Tolimieri/Orr and Daubechies et al. \cite{Jan95,DLL95,TO95}. The
main result of Wexler/Raz consists in our setting of the following condition:
\par
Let $\Lambda^0$ be a lattice in $\R^{2d}$. A pair $(g,\gamma)\in
M^{p.q}_m\times M^{p',q'}_{1/m}(\R^d)$ satisfies the {\it Wexler-Raz
condition} with respect to $\Lambda^0$, if
\begin{equation}
   |\Lambda|^{-1}\langle \gamma,\pi(\lambda^0)g\rangle=\delta_{0,\lambda^0},
\end{equation}
where $\delta_{0,\lambda^0}$ denotes the Kronecker delta for the set
$\Lambda^0$. In terms of Gabor systems the Wexler-Raz condition expresses the
{\it biorthogonality} of the two sets ${\mathcal G}(g,\Lambda^0)$ and
${\mathcal G}(\gamma,\Lambda^0)$ to each other on $L^2(\R^d)$.
\par
The importance of the Wexler-Raz condition arises from the fact, that under
certain assumptions it characterizes all dual atoms of a given Gabor frame
${\mathcal G}(g,\Lambda)$.
\par
The following theorem is the proposed extension of Feichtinger and
Zimmermann's result \cite{FZ98}.
\begin{thm}
 Let $\Lambda$ be a lattice in $\R^{2d}$ and let $(g,\gamma)$ be a dual pair in $M^{p.q}_m\times M^{p',q'}_{1/m}(\R^d)$. Then the following holds:
 \begin{enumerate}
   \item[(1)] (Wexler-Raz Identity)
      \begin{equation}
         S_{g,\gamma,\Lambda}f=|\Lambda|^{-1}S_{f,\gamma,\Lambda^0}g~~~\text{in}~~~M^{\infty}(\R^d)
      \end{equation}
      for all $f\in M^{1}(\R^d)$.
   \item[(2)] (Janssen Representation)
     \begin{equation}
         S_{g,\gamma,\Lambda}=|\Lambda|^{-1}\sum_{\lambda^0\in\Lambda^0}V_{\gamma}g(\lambda^0)\pi(\lambda^0)
      \end{equation}
      is a bounded operator from $M^{1}(\R^d)$ to $M^{\infty}(\R^d)$ and the series converges unconditionally in the strong sense.

 \end{enumerate}
\end{thm}
The proof is just a reformulation of the FIGA and the arguments of
\cite{FZ98} are also valid in our situation.
\par
The usefulness of weakly dual pairs relies on the fact, that it is
equivalent to the Wexler-Raz condition.
\begin{thm}[Feichtinger-Zimmermann]

Let $\Lambda$ be a lattice in $\R^{2d}$. Then a pair $(g,\gamma)$ in
$M^{p.q}_m\times M^{p',q'}_{1/m}(\R^d)$ is weakly dual with respect to
$\Lambda$ if and only if $(g,\gamma)$ satisfies the Wexler-Raz condition with
respect to $\Lambda^0$.
\end{thm}
The proof of Feichtinger/Zimmermann can again be adapted to our situation
\cite{FZ98}.

\par
{\bf Acknowledgement:} The authors want to thank W. Czaja and H. Rauhut  for many useful
suggestions. The second named author wants to thank the Max Planck Institute
of Mathematics at Bonn for their hospitality and Prof. Y. Manin for his kind
invitation, because large part of the manuscript was written during this
stay.

\flushleft{
}

\end{document}